\documentstyle[12pt,amssymb,epsfig]{amsart}
\newcommand{\be}{\begin{enumerate}}
\newcommand{\ee}{\end{enumerate}}
\newcommand{\bd}{\begin{description}}
\newcommand{\ed}{\end{description}}
\newcommand{\bi}{\begin{itemize}}
\newcommand{\ei}{\end{itemize}}
\begin{document}

\newtheorem{theorem}{Theorem}[section]
\newtheorem{proposition}[theorem]{Proposition}
\newtheorem{lemma}[theorem]{Lemma}
\newtheorem{corollary}[theorem]{Corollary}
\newtheorem{definition}[theorem]{Definition}
\newtheorem{remark}[theorem]{Remark}
\newtheorem{conjecture}[theorem]{Conjecture}
\newcommand{\eproof}{\rule{0,2cm}{0,2cm}}
\newcommand{\nat}{I\!\!N}

\begin{center}{\bf On $q$-Gaussians and Exchangeability}

\bigskip
Marjorie G. Hahn$^1$, Xinxin Jiang$^2$, and Sabir Umarov$^1$

 \medskip
 $^1$Department of Mathematics, Tufts University, Medford, MA 02155 USA
 $^2$Department of Mathematics and Computer Science, Suffolk University, Boston, MA 02114  USA

\bigskip

\textbf{Abstract}
\end{center}

 {\small The $q$-Gaussians are discussed from the point of view of
variance mixtures of normals and exchangeability.  For each
$-\infty<q< 3$, there is a $q$-Gaussian distribution that maximizes
the Tsallis entropy under suitable constraints. This paper shows
that $q$-Gaussian random variables can be represented as variance
mixtures of normals. These variance mixtures of normals are the
attractors in central limit theorems for sequences of exchangeable
random variables; thereby, providing a possible model that has been
extensively studied in probability theory. The formulation provided
has the additional advantage of yielding process versions which are
naturally  $q$-Brownian motions. Explicit mixing distributions for
$q$-Gaussians should facilitate applications to areas such as option
pricing.  The model might provide insight into the study of
superstatistics. }


\vspace{1.5cm}

\section{Introduction}

Developments in nonextensive statistical mechanics based on an
entropy proposed by Constantino Tsallis (1988) gave rise to the
$q$-Gaussian distributions which are being applied in numerous
research areas,  see Gell-Mann and Tsallis (2004), Boon and Tsallis
{\sl at al} (2005), and
Tsallis (2009).
Applications are based on experimental, computational, and
analytical results.
By definition, for $-\infty <q<3$, the
$q$-Gaussian density takes the form $g_q(x)=C_q
[1-(1-q)x^2]^{1/(1-q)}$ where $C_q$ is specified in Section 2.
Theoretical under-pinnings of the $q$-Gaussians are founded  on a
novel $q$-algebra.
The purpose of this paper is to propose a stochastic model within
the framework of usual algebra that is consistent with many of the
phenomena people are attempting to model by the $q$-Gaussian
distributions with tails that decay as a power of $x$ (those with
$1< q<3$) as well as the Gaussian distribution ($q=1$). The model
not only facilitates theoretical investigations, but helps to
clarify why $q$-Gaussian distributions arise in certain physical
applications.

For each $-\infty<q<3$, Tsallis (1988) defined a  $q$-Gaussian
distribution  to be that distribution which, under certain
conditions, maximizes the Tsallis $q$-entropy:  $$\displaystyle
S_q=k{1-\sum_{i=1}^W p_i^q\over q-1}, $$  where $W\in \mathbb{N}$ is
the total number of possible micropscopic configurations, $\{p_i\}$
are the associated configuration probabilities with
$\sum_{i=1}^Wp_i=1$, and $k$ is a conventional positive constant.
When $q=1$ the definition is understood via the limit, in which case
the $1$-entropy recovers the Boltzmann-Gibbs entropy. The
$1$-Gaussian is simply the usual Gaussian distribution. The
$q$-Gaussian distributions have tails that decay on the order of
$|x|^{-2/(q-1)}$ when $1<q <3$ and have bounded support when $q<1$.

To a large extent, prevalence of the Gaussian distributions in
theory and applications stems from their roles as attractors, e.g.
via the classical central limit theorem for sequences of
independent, identically distributed random variables with finite
second moments. The dependence structure that should lead to
$q$-Gaussian limits has remained elusive. The phrase ``globally
correlated" is often attached to phenomena in statistical mechanics
to which Tsallis entropy is applied.

There have been efforts to
understand nonextensivity and the role of $q$-Gaussians from a
probability point of view. For example, Umarov {\sl et al} (2008)
proves $q$-central limit theorems using $q$-algebra and a notion of
$q$-independence.
Since
the system developed is surprisingly nice mathematically, it is
interesting
to understand in the sense of usual probability the intuitive
meanings of $q$-algebra and $q$-independence.
Using usual algebra and usual notions of dependence, our paper establishes that  $q$-Gaussians, with $1\le
q<3$, have a role as natural attractors since they are a subset of the possible limits in a general central limit theorem for dependent random variables. This provides an explanation for the wide occurrence of at least these $q$-Gaussians.

Our initial investigations were stimulated by the common features of two specific mathematical models for which it is possible to analytically verify whether or not the weak limit is $q$-Gaussian.
Vignat and Plastino (2007) suggests a model
which combines
randomness in the normalizer with each i.i.d. random variable, and
shows that the weak limit is $q$-Gaussian. Their model can be viewed
differently, namely as a special case of a sequence of exchangeable
random variables.
Marsh {\sl et al} (2006) provides practical models by applying
Leibnitz triangles, in which the triangular arrays of random
variables are rowwise exchangeable, however
weak limits are not $q$-Gaussian; see Hilhorst and Schehr (2007).
Both examples are discussed further in Section 2 since they
motivated our consideration of exchangeability as a possible
probability model that might be connected to some applications being
considered using Tsallis nonextensivity theory.

From a different but related direction, this paper is also
stimulated by a possible connection between the proposed model and
superstatistics.  Beck and Cohen (2003) introduces the concept of
superstatistics for generalized Boltzman factors derived from
systems that evolve in complex environments.  The expression for a
superstatistic can be interpreted as  a variance mixture of normals
(VMON) (or scale mixture of normals) after standardization.
Moreover, in a specific case, the superstatistic reduces to the
Tsallis statistic (Beck (2001)). VMONs are the attractors in a
central limit theorem for exchangeable sequences, see Jiang and Hahn
(2003). Furthermore, the reason behind studying superstatistics,
such as a non-homogeneous background where a mechanical system
expands, strongly suggests exchangeability in modeling the
underlying system. De Finetti's theorem ( Chow and Teicher (1978))
characterizes an infinite sequence of exchangeable random variables
as a mixture of i.i.d. sequences, i.e. an infinite exchangeable
sequence is conditionally an i.i.d. sequence. Thus, if a mechanical
system evolves in a complex, non-homogeneous background with the
usual central limit theorem (for i.i.d. sequences) holding in each
homogeneous part of the background, then the limiting process has
marginals which are mixtures of normals. Further discussion of
superstatistics appears at the end of Section 2.

Similar arguments can be made in finance.  Viewing market movements
as responses to different specific information (interest rate, macro
financial data, earnings, etc.), then conditioning on each specific
type of information, yields a market which follows a random walk.
Consequently, distributions that can model asset returns could
naturally be mixtures of normals. Empirical evidence (Hall {\sl et
al} (1989), Gribbin {\sl et al} (1992), and Kon (1984))
 has shown that discrete VMONs fit asset returns data better
than stable distributions or than the Student model.  Stables have
been studied more extensively than discrete VMONs.  Moreover,
research literature on applications of continuous VMONs seems to be
quite limited.

A drawback to applications of VMONs has always been that the mixing
distributions are not generally easy to obtain from data.
Expectation Maximizing (EM) algorithms are usually employed.
However, when the number of Gaussians in the mixture is large, or in
many cases infinite, computations become very complicated. This
paper shows that $q$-Gaussians are VMONs when $1\le q< 3$ (and not
when $q<1$). Furthermore, the mixing distributions are calculated
explicitly; thereby, facilitating application of these models. For
example, the option price obtained from a VMON is the mixture of the
option prices when the variance is fixed.  Therefore, if one
believes that $q$-Gaussians are good models for financial returns,
then their option prices can be calculated simply by mixing the
prices obtained from the classical Black-Scholes model with the
known mixing distributions. Furthermore, studying mixing
distributions might provide insight into the ``background". For
instance, a belief that information drives stock prices suggests
viewing the mixing distributions as models for the impact of
information.

The examples discussed above are each special cases of the model
proposed in this paper.  Our model of $q$-Gaussian distributions
based on exchangeability, which is restricted to $1\le q<3$, is
specified using {\it usual algebra},  rather than $q$-algebra, and
has the following features:
\begin{itemize}
\item[1)] specification of  $q$-Gaussian random variables as specific variance mixtures of normals;
\item[2)] associated central limit theorems for dependent exchangeable random variables with the $q$-Gaussians as natural attractors;
\item[3)] associated stochastic processes which in a natural sense are $q$-Brownian motions.
\end{itemize}
Within this model, as the latter two features indicate, the notion
of ``global correlation" is specified as exchangeability.  From this
point onward we consider the appearance of q-Gaussians only as
variance mixture of normals under the exchangeability concept, not
discussing their appearance from other concepts.

\subsection{Organization of the paper}
 Section 2 makes the connection
between $q$-Gaussians and variance mixtures of normals as well as
obtains a complete specification of the mixing distributions, which
facilitates applications.  Superstatistics are discussed as an
example.  Section 3 provides a detailed mathematical description of
exchangeability and establishes central limit theorems for
exchangeable sequences and rowwise exchangeable triangular arrays.
The examples motivating our consideration of exchangeability are
discussed further in this section. Our concept of $q$-Brownian
motion is introduced in section 4 followed by some comparisons with
a stochastic process in Borland (1998) that has $q$-Gaussian
marginal distributions. Section 5 is the conclusion.

\subsection{Conventions}  Throughout the paper, $Z$ will be used
specifically for a standard normal random variable that is
independent of any other random variables. $\to^{\cal L}$, $\to^P$ ,
$\to^{a.s.}$ stand respectively for convergence in distribution, in
probability, and almost surely.  It is worth re-emphasizing that all
further discussion relies on {usual algebra} rather than
$q$-algebra..
\medskip

\section{$q$-Gaussians and variance mixtures of normals}

The {\it $q$-Gaussian distributions} form a one-parameter family of
distributions for $-\infty<q<3$ with densities specified by
$g_q(x)=C_q [1-(1-q)x^2]^{1/(1-q)}, $ where the normalizing constant
$C_q$ is given by
\begin{equation*}
C_{q}^{-1}=\left\{
\begin{array}{ll}
{\ \frac{2}{\sqrt{1-q}}\int_{0}^{\pi
/2}(\cos\,t)^{\frac{3-q}{1-q}}dt}=\frac{2 \sqrt{\pi }\,\Gamma
\bigl({\frac{1}{{1-q}}}\bigr)}{(3-q)\sqrt{1-q}\,\Gamma
\bigl({\frac{{3-q}}{{2(1-q)}}}\bigr)}, & -\infty <q<1, \\
\sqrt{\pi }, & q=1, \\
\frac{2}{\sqrt{q-1}}\int_{0}^{\infty
}(1+y^{2})^{\frac{{-1}}{{q-1}}}dy=\frac{ \sqrt{\pi }\,\Gamma
\bigl(\frac{3-q}{2(q-1)}\bigr)}{\sqrt{q-1}\,\Gamma \bigl(
{\frac{1}{{q-1}}}\bigr)}, & 1<q<3\,. \\
&
\end{array}
\right.
\end{equation*}
When $q=1$, the expression for $g_q$ is understood by taking limits
and yields the  standard normal density function.
 The support of the $q$-Gaussian is $(-\infty, \infty)$ if $1\le q< 3$ and the compact set  $\displaystyle{\left[-\frac{1}{\sqrt{1-q}},\frac{1}{\sqrt{1-q}}\right]}$, if $-\infty<q<1$.

A {\it variance or (scale) mixture of normal distributions (VMON)}
by definition has characteristic function of the form
$\phi(t)=\int_0^{\infty} \exp(-t^2u/2)d H(u),$ with $H$ a
distribution function on $[0,\infty)$. The corresponding density is
$$f(x)=\int_0^{\infty} (2\pi u)^{-1/2}\exp(-x^2/(2u))d H(u).$$ Each
VMON has a representation of the form $VZ$ where $V>0$ a.s., and $V$
is independent of $Z$. The VMONs include many commonly used
distributions such as the symmetric stable distributions, the
Cauchy, Laplace, double exponential, logistic, hyperbolic, and
Student distributions and their mixtures, plus many others. See
Keilson and Steutel (1974) and Gneiting (1997).

The connection between the $q$-Gaussians and the VMON can be made by
applying the theory of completely monotone functions. A function $h$
is {\it completely monotone
} on $(0,\infty)$ if and only if $(-1)^nh^{(n)}(x)\geq 0$ for $x>0$,
and $n=0, 1, 2, ...$

Andrews and Mallows (1974) shows that a symmetric density function
$f_X(x)$ is a variance mixture of normals if and only if
$f_X(\sqrt{x})$ is completely monotone.  It is easy to show that
$g_q(\sqrt{x})$ is completely monotone for $q>1$ and not completely
monotone for $q<1$, which leads to the following statement.

\begin{theorem} $q$-Gaussians are variance mixtures of normals when $1<q<3$
and not mixtures of normals when $q<1$.
\end{theorem}

Since not every VMON is a $q$-Gaussian, it remains to identify the
mixing distributions that yield the $q$-Gaussians for $1<q<3$.  In
Beck (2001) the $q$-Gaussians are identified as $\frac{1}{a} Z$
where $a^2$ has a $\chi^2$ distribution with $q$ derived from the
degrees of freedom. Below we provide a direct proof using the
Laplace transform technique when the inverse Laplace transform of
the density is known.   It illustrates a method  applicable to
identification of unknown mixing measures from other variance
mixtures of normals as well.

\begin{theorem} \label{gGaussThm}   The $q$-Gaussian density
for $1<q<3$, can be expressed as the following variance mixture of
normal densities:
$$g_q(x)\equiv f_{V\cdot Z, q}=\int_0^\infty {1\over\sqrt{2\pi}v}\exp\left({-x^2\over 2v^2}\right)
f_V(v)\, dv$$ where the mixing measure $$dH(v)=f_V(v)\,
dv=C_{V,q}\exp\left(-{1\over 2(q-1)v^2}\right) v^{-{2\over q-1}}\,
dv,$$ with \quad $\displaystyle{C_{V,q}^{-1}=\Gamma\left({3-q\over
2(q-1)}\right)\cdot {1\over 2}\cdot [2(q-1)]^{3-q\over 2(q-1)}}$.
\end{theorem}

{\it Proof.} Let $\mathnormal {L}$ denote the Laplace transform.
Since $\displaystyle{{\mathnormal L}\left(\exp(-bt)\cdot
t^\alpha\right)={\Gamma(\alpha+1)\over (\xi+b)^{\alpha+1}}}$ for $
\alpha>-1$,  the choice of $\displaystyle{\alpha={1\over
q-1}-1={2-q\over q-1},\quad b={1\over q-1}}$ yields
$${\mathnormal L}\left(\exp\left(-{t\over q-1}\right)\cdot
t^{2-q\over q-1}\right)={\Gamma\left({1\over q-1}\right)\over
\left(\xi+{1\over q-1}\right)^{1\over q-1}}={\Gamma\left({1\over
q-1}\right)(q-1)^{1\over q-1}\over [1+(q-1)\xi]^{1\over q-1}}.$$
Equivalently,
$$
\displaystyle{C'_q [1-(1-q)x^2]^{1/(1-q)}=\int_0^\infty
\exp(-x^2t)\, d H(t)},
$$
where  $\displaystyle{dH(t)=\exp\left(-{t\over q-1}\right)\cdot
t^{2-q\over q-1}\, dt \quad {\rm and}\quad C'_q= \Gamma\left({1\over
q-1}\right)(q-1)^{1\over 1-q}}$.
The following routine calculations verify the claim where the
following substitutions are used: $v^2\to v$ in the second equality,
$1/v\to u$ in the third equality, and $ \left({1\over
2}\left(x^2+{1\over q-1}\right)u\to v\right)$ in the fourth
equality:

\begin{eqnarray}
\nonumber &&
\int_0^\infty {1\over\sqrt{2\pi}v}\exp\left({-x^2\over 2v^2}\right)
f_V(v)\, dv  \\
 \nonumber&&=C_{V,q}\int_0^\infty {1\over\sqrt{2\pi}v}\exp\left({-x^2\over 2v^2}\right)
  \exp\left(-{1\over 2(q-1)v^2}\right) v^{-{2\over q-1}}\,
dv \\
\nonumber&&= {1\over\sqrt{2\pi}}C_{V,q}{1\over 2}\int_0^\infty
\exp\left(-{1\over 2}\left(x^2+{1\over q-1}\right)
\cdot {1\over v}\right)\left({1\over v}\right)^{q\over q-1}dv\\
\nonumber&&= {1\over\sqrt{2\pi}}C_{V,q}{1\over 2}\int_0^\infty
\exp\left(-{1\over 2}\left(x^2+{1\over q-1}\right)
u\right)u^{{1\over q-1}-1}du \quad\\
 \nonumber&&={1\over\sqrt{2\pi}}C_{V,q}{1\over 2}\int_0^\infty \exp(-v)\cdot v^{{1\over
 q-1}-1}\cdot \left({1\over 2}\left(x^2+{1\over
 q-1}\right)\right)^{-1\over q-1} dv 
 \end{eqnarray}

\begin{eqnarray}
\nonumber&&={1\over\sqrt{2\pi}}C_{V,q}{1\over 2}\Gamma\left({1\over
q-1}\right)\left({1\over 2}\left(x^2+{1\over
q-1}\right)\right)^{-1\over q-1} \\
\nonumber&&= {1\over\sqrt{2\pi}}C_{V,q}{1\over 2}2^{1\over
q-1}\Gamma\left({1\over q-1}\right)(q-1)^{1\over
q-1}\left(1+(q-1)x^2\right)^{-1\over q-1}.
\end{eqnarray}
It is not hard to verify that $${1\over\sqrt{2\pi}}C_{V,q}{1\over
2}2^{1\over q-1}\Gamma\left({1\over q-1}\right)(q-1)^{1\over q-1}=
C_q,$$ which completes the proof. \eproof

In Theorem \ref{gGaussThm}, $V$ is 1 over the square root of a
$\chi^2$ distribution with the number of degrees of freedom being
$\frac{2}{q-1} -1$. Also notice that when $q=2$, the $q$-Gaussian is
Cauchy and
$$f_V(v)={\sqrt 2\over\sqrt \pi}{1\over v^2}\exp\left({-1\over
2v^2}\right).$$

When $q<4/3$, $EX^2<\infty$. This situation is important in finance
or risk management where the variance is often used  to quantify the
risk. In that situation, the $q$-Gaussians are superior to
non-Gaussian stable distributions which fail to have finite second
moments. When $q<3/2$, $EX<\infty$. This situation is important in
finance or insurance when measuring mean returns or mean expenses.

The mixing distributions are a type of generalized inverse Gaussian
distribution, and the $q$-Gaussians are a type of generalized
hyperbolic distribution. It is worth noting that both
 generalized inverse Gaussian and generalized hyperbolic distributions
  are infinitely divisible.
 Since infinitely divisible laws are the only attractors for triangular
 arrays of independent random variables, this fact might help explain why
 $q$-Gaussians should be often observed without self-normalization (i.e. without pre-processesing).

\subsection{Example: Superstatistics}

As a test particle moves from cell to cell, its velocity $v$
satisfies a Langevin equation $$\dot{v}=-\gamma v+\sigma \dot
B_t,\eqno(1)$$ where $B_t$ is a standard Brownian motion. Instead of
$\gamma$ and $\sigma$ being deterministic, as in classical
statistical physics, Beck and Cohen (2003) let
$\beta=\gamma/\sigma^2$ be a random variable with density
$f(\beta)$. In this setting, the widely used classical Boltzman
factor $e^{-\beta E}$ takes a generalized form called a
superstatistic:
$$B(E)=\int_0^\infty k(\beta)e^{-\beta E}d\beta, \eqno(2)$$
where $E$ represents the energy of a microstate associated with each
cell. Notice that the superstatistic in  (2) is exactly a variance
mixture of normals after standardization. The special case where
$k(\beta)$ (the density of $1/V^2$ in our Theorem \ref{gGaussThm})
is the density of a  $\chi^2$ random variable  yields the Tsallis
statistic (Beck (2001)).

\section{ Central limit theorems for exchangeable random variables
with q-Gaussian limits}

A brief introduction to exchangeability is required. For details,
see Chow and Teicher (1978).

A sequence of random variables $\{X_n,\ n\geq 1\}$ defined on some
probability space is said to be exchangeable if for each $n$,
$$P(X_1\leq x_1,...,X_n\leq x_n)=P(X_{\pi(1)}\leq
x_1,...,X_{\pi(n)}\leq x_n)$$ for any permutation $\pi$ of
$\{1,2,...,n\}$ and any $x_i\in {\bf R}$, $i=1,...,n.$

By de Finetti's Theorem, an infinite sequence of exchangeable random
variables is conditionally i.i.d.$\!$ given the $\sigma$-field
${\cal G}$ of permutable events. Furthermore,
there exists a regular
conditional distribution $P^\omega$ for $X_i$ given ${\cal G}$ such
that for each $\omega\in\Omega$ the coordinate random variables
$\{\xi_n\equiv\xi^\omega_n, n\geq 1\}$ (called mixands) are i.i.d.
Hence, for each natural number $n$, any Borel function $f:{\bf
R}^n\to {\bf R}$, and any Borel subset of ${\bf R}$,
\[
P(f(X_1,\ldots , X_n)\in B) =\  \int_{\Omega} P(f(X_1,\ldots ,
X_n)\in B\vert {\cal G}) dP \]
\[
\hspace{4.2cm} = \ \int_{\Omega} P^{\omega}(f(\xi_1,\ldots ,
\xi_n)\in B) dP.
\]
The mixands are allowed to be on a different probability space than
the $X_i$'s.

For exchangeable sequences, the dependence never dies in contrast to
weakly dependent sequences. Hence if the covariance exists, it does
not change along the sequence. Moreover, for an infinite
exchangeable sequence the covariance is always non-negative.

\subsection{Example.} The model proposed in Marsh {\sl et al}
(2006), considers $N$ identical and distinguishable, but not
necessarily independent binary subsystems. Let $r_{N,n}$ be the
probability that there are $n$ subsystems in state 1, which is given
by the Leibnitz rule:
$$r_{N,
n}+r_{N, n+1}=r_{N-1, n}.$$ Since construction of the model only
considers the number of subsystems in state 1, the order of 1's and
0's does not matter, a typical property of exchangeable sequences.
However, it is unclear whether each row of $N$ 0-1 random variables
can be embedded in an infinite sequence of exchangeable random
variables.
Rodriguez {\sl et al} (2008) and
Hanel {\sl et al} (2009) construct exchangeable models that can be
embedded into an infinite sequence of exchangeable random variables,
again using the Leibnitz rule.

Turning to central limit theorems, the following simplified and
adapted version of Theorem 2.1 in Jiang and Hahn (2003) suffices for
the needs of this paper. The proof is provided for completeness.

\begin{theorem} \label{first}Let $\{ X_{n},\ n \geq 1 \}$ be an infinite sequence of
exchangeable random variables where $X_1$ has mean 0 if the mean
exists or is symmetric otherwise. Assume $0<E\,\xi_1^2<\infty$ a.s.
Then either
$$i)\quad {\sum_{i=1}^n X_i\over \sqrt{n}} \to^{\cal L} V_1\cdot Z$$
or $$ii)\quad {\sum_{i=1}^n X_i\over n}\to^{a.s.} V_2,$$ where
$V_1=\sqrt{Var(\xi_1)}=\sqrt{Var(X_1|\,{\cal G})}$ and $V_2=E\,
\xi_1=E(X_1|\,{\cal G})$.
\end{theorem}

{\bf Proof.} By de Finetti's Theorem, for any real $x$,
$$P\left({\sum_{i=1}^n X_i\over \sqrt n}\leq
x\right)=\int_\Omega P^\omega\left({\sum_{i=1}^n\xi_i-nE\xi_1\over
\sqrt n \cdot \sigma(\xi_1)}\sigma(\xi_1)+\sqrt nE\xi_1\leq
x\right)dP,$$ where $\sigma(\xi_1)$ denotes the standard deviation
of $\xi_1$. If $E\xi_1=0$ a.s., then the classical central limit
theorem holds for the mixands and case i) of the theorem holds. If
$E\xi_1$ is not almost surely 0, replace
 the $\sqrt n$ in the previous equation by $n$. Then the first summand
 inside the regular conditional probability on the right side integral
 converges to zero almost surely and case ii) of the theorem holds.
\eproof

\begin{remark} The condition put on the mixands $0<E\,\xi_1^2<\infty$
a.s. is generally weaker than assuming $X_i$'s have finite second
moments.
\end{remark}

\begin{remark} Case i) of Theorem \ref{first} shows that $q$-Gaussians
when $1<q<3$ are among the possible limits in a central limit
theorem for exchangeable sequences. However, they are not the only
ones.  It seems more natural to consider $q$-Gaussians as limit
distributions in case i), since theoretically the range of
distributions in case ii) is vast. For example, let $\epsilon_i$'s
be i.i.d. with mean 0 and standard deviation 1, and $Y$ be any
random variable that is independent of all $\epsilon_i$'s. Then
$\{Y+\epsilon_i, i\geq 1\}$ is a sequence of exchangeable random
variables and ${\sum(Y+\epsilon_i)/n}\to^{a.s.} Y$.
\end{remark}

\begin{remark}  In case i), $EX_1X_2=0$ if it exists, which
means that $q$-Gaussians can be attractors of uncorrelated but not
independent exchangeable sequences of random variables. In the
second case, $EX_1X_2\not=0$ if it exists and in this case the
central limit theorem is really the strong law of large numbers for exchangeable random variables.
This case explains why the examples in Rodriguez {\sl et al} (2008) and Hanel {\it et al} (2009) achieve $q$-Gaussian limits with normalizers at the rate of $n$.
\end{remark}

\begin{remark} Theorem 4.1 of Vignat and Plastino
(2007)
says
that if $Y_i$  are i.i.d. in the domain of normal attraction of
$Z$,and $a^2$ has a $\chi^2$-distribution (correcting a typo), then
${1\over a\sqrt{n}}\sum_1^n Y_i$ converges in distribution to a
$q$-Gaussian distribution. Their theorem in one dimension is a
special case of Theorem \ref{first} in this paper with $X_i=V_1\cdot
Y_i$ where $V_1=\frac{1}{a}$. However, in general, a sequence of
exchangeable random variables does not necessarily have the form
$V\cdot Y_i$.
\end{remark}

The next theorem is a triangular array version of Theorem
\ref{first}, which seems to be new.

\begin{theorem} \label{second} Let $\{X_{n,i}, i\leq n\leq n, n=1, 2, \dots\}$ be
a triangular array of rowwise exchangeable random variables that can
be embedded into infinite exchangeable sequences. Assume $X_{n,i}$'s
are centered or symmetric with $0<E\,\xi_{n,1}^2<\infty$ in
probability for each $n$, and $Var(\xi_{n,1})\to V_1^2$ in
probability when $n\to\infty$, with $V_1>0$ a.s. Then
$${\sum_{i=1}^n X_{n,i}\over\sqrt {n}}\to^{\cal L}V_1\cdot Z+V_2$$
when $\sqrt nE\xi_{n,1}\to V_2$, and
$${\sum_{i=1}^n X_{n,i}\over n^\alpha}\to^P V_3$$  with $\alpha>1/2$,
when $n^{1-\alpha} E\xi_{n,i}\to^P V_3$  as $n\to\infty$.
\end{theorem}

{\bf Proof.} By de Finetti's Theorem, for any real $x$,
$$P\left({\sum_{i=1}^n X_{n,i}\over \sqrt n}\leq
x\right)=\int_\Omega
P^\omega\left({\sum_{i=1}^n\xi_{n,i}-nE\xi_{n,1}\over \sqrt n \cdot
\sigma(\xi_{n,1})}\sigma(\xi_{n,1})+\sqrt nE\xi_{n,1}\leq
x\right)dP.$$ Using the definition of weak convergence in
probability (Hahn and Zhang (1998)), and the conditions that $\sqrt
nE\xi_{n,1}\to^P V_2$,
 $\sigma(\xi_{n, 1})\to^P V_1$ with $V_1>0$
a.s., and that the mixands have finite second moments almost surely,
yields, when $n\to \infty$
$$P\left({\sum_{i=1}^n X_{n,i}\over \sqrt n}\leq
x\right)\to \int_\Omega P^\omega( V_1Z+V_2\leq x)dP=P(V_1Z+V_2\leq
x).$$

When $n^{1-\alpha} E\xi_{n,1}\to^P V_3$ with $\alpha>1/2$ and
$V_2\not = 0 $ a.s., again, by de Finetti's Theorem,
$$P\left({\sum X_{n,i}\over n^\alpha}\leq x\right )=\int_\Omega
P^\omega\left({\sum_{i=1}^n\xi_{n,i}-nE\xi_{n,1}\over n^\alpha \cdot
\sigma(\xi_{n,1})}\sigma(\xi_{n,1})+n^{1-\alpha} E\xi_{n,1}\leq
x\right)dP.$$ Since $${\sum_{i=1}^n\xi_{n, i}-nE\xi_{n,1}\over
n^\alpha\, \sigma(\xi_{n,1})}\to^P 0$$ when $\alpha>1/2$ and
$\sigma(\xi_{n,1})\to^P V_2$,
$${\sum_{i=1}^n\xi_{n,i}-nE\xi_{n,1}\over n^\alpha \cdot
\sigma(\xi_{n,1})}\sigma(\xi_{n,1})\to^P 0,$$ and the second part of
the theorem is proved. \eproof

\begin{remark}  When $X_{n,i}$'s are identically distributed, then
$V_2=0$ and $\alpha=1$, which is consistent with Theorem
\ref{first}.
\end{remark}

\section{ $q$-Brownian Motion}

The literature on $q$-Gaussians thus far has failed to allow the
construction of a $q$-analogue of Brownian motion.  Our
representation of $q$-Gaussians for $q>1$ in terms of variance
mixtures of normals, using usual probabilistic and algebraic
notions, leads naturally to the definition of a  process that might
naturally be called $q$-Brownian motion.  However, we first require
some definitions.

\begin{definition}
A stochastic process is called {\it exchangeable} if it is
continuous in probability with $X_0=0$ and such that the increments
over disjoint intervals of equal length form an exchangeable
sequence.
\end{definition}
\begin{definition}
$X$ has {\it conditionally independent, stationary increments, given
some $\sigma$-field}, if \underbar{both} properties of the
increments are conditionally valid for any finite collection of
disjoint intervals of the same length.
\end{definition}
These two definitions are connected by a theorem of B\"uhlmann which
characterizes exchangeable processes on $R^+$.   [See e.g. Theorem
9.21 of O. Kallenberg (2002).]

\begin{theorem} Let the process $(X_t)_{t\ge 0}$  be
$R^d$-valued and continuous in probability with $X_0=0$.  Then $X$
is exchangeable if and only if it has conditionally independent,
stationary increments given some $\sigma$-field.
\end{theorem}

Using these definitions, we can now define a $q$-Brownian motion.

\begin{definition}   A {\it $q$-Brownian motion}, $B^q_t$, with
$1\le q<3$ is a stochastic process having the following properties:

1) conditionally independent, stationary increments;

2) all increments are $q$-Gaussian;

3)  a.s. continuous sample paths.
\end{definition}

\begin{theorem} (Existence)  For each $1\le q<3$,  a
$q$-Brownian motion exists.
\end{theorem}

{\it Proof.}  Let $B_t$ be a standard Brownian motion.  When $q=1$,
$B^1_t=B_t$ is  a process satisfying all three conditions with the
conditional $\sigma$-field being the trivial $\sigma$-field. For
$1<q<3$,  let $V_q$ be a random variable with the density $f_q$
given in Theorem \ref{gGaussThm}.  Then $B^q_t=V_q \cdot B_t$ has
conditionally stationary and independent increments given the value
of $V_q$.  All increments are $q$-Gaussian.  Furthermore, a mixture
of processes with continuous sample paths, has continuous sample
paths. \eproof

\subsection{Example. Comparison with the Borland process}

Borland (1998) provides another process with  $q$-Gaussian
marginals. Let $Y_t$ be the log returns of stock prices that follow
the stochastic differential equation $dY_t=\mu dt+\sigma d\Omega_t$,
where $\Omega_t$ evolves according to
$d\Omega_t=P(\Omega_t)^{(1-q)/2}d B_t$. The evolution of the
probability distribution $P$ of $Y_t$ is nonlinear according to the
nonlinear Fokker-Planck (forward Kolmogorov) equation
$${\partial\over \partial t}P(x, t|y)={1\over
2}{\partial^2\over \partial x^2}\left[P^{2-q}(x, t|y)\right].$$
Solutions for $P(x, t|\, y)$ are given by $q$-Gaussians for each
fixed $t$, as established in Borland (1998) or using the $q$-Fourier
transform as in Umarov and Queiros (2009).

Our $q$-Brownian motion has conditionally, stationary and
independent increments, each of which is $q$-Gaussian. This is
expressed in the fact that its transition density is expressed in
the form
$$P(t, x|\, y)=\int_0^{\infty}p(t, x|\, y,\, v) f_V(v) dv,$$
where $p(t, x| \, y,  v)$ for each fixed $v$ satisfies the linear
Fokker-Planck equation
$${d\over dt}p
(t, x|\, y, v)={1\over 2}v^2{\partial^2 p
(t, x|\, y, v)\over \partial x^2},$$ and the occurrence of the
transition probability for a particular $v$ is weighted according to
the density $f_V(v)$ specified in Section 2.

The Borland process clearly differs from our $q$-Brownian motion.
Even though it has stationary increments (Umarov and Queiros (2009), its increments are not
conditionally independent. It does have continuous paths and
$q$-Gaussian marginals.  It is important to further clarify the nature of the increments for the
Borland process for financial and other applications.

\section{Conclusion}

This paper provides a probabilistic model for $q$-Gaussian
distributions with $1\le q <3$ based on exchangeability as one of
the possible notions of ``global correlation.''  The model is
consistent with many of  the phenomena for which $q$-Gaussian
distributions are being used based on empirical evidence.  In
particular, $q$-Gaussian distributions are variance mixtures of
normals when $1\le q <3$ and not when $q<1$. Explicit mixing
distributions are provided, which should extend further application
of $q$-Gaussian distributions.
An explanation for the wide occurrence of these $q$-Gaussians is their role as attractors via central limit theorems
for exchangeable sequences and triangular arrays.
A natural $q$-Brownian motion is
defined which can be viewed as an alternative driving process to the
Borland (1998) process.  The increments of the two processes have
different characteristics and thus model different phenomena.  The
paper also makes a connection with superstatistics which are
variance mixtures of normals after normalization. The Langevin
equations that yield superstatistics can be viewed as stochastic
differential equations driven by exchangeable processes.

\vspace{1cm}

\noindent {\bf References}

 \bigskip \frenchspacing

\be

 \item{} Andrew D. F. and Mallows C. L.  (1974)
{\it Scale mixtures of  normal distributions.} {\it J. Royal Stat.
Soc. Series B.} {\bf 36}, 99-102.

\item{} Beck C. (2001)
{\it Dynamical foundations of nonextensive statistical mechanics.}
{\it Phys. Rev. Lett.} {\bf 87}, 180601.

 \item{} Beck C. and Cohen E. G. D. (2003)
 {\it Superstatistics.} {\it  Physica A.} {\bf 322}, 267-275.

 \item{} Boon J. P. and Tsallis C. (2005)
 {\it Nonextensive Statistical Mechanics: New Trends, New Perspectives,
  Europhysics News,} {\bf 36}, 6, European Physical Society.

\item{} Borland L. (1998)
{\it Microscopic dynamics of the nonlinear Fokker-Plank equations: A
phenomenological model.} {\it Phys. Rev. E} {\bf 57,} 6634-6642.

\item{} Chow Y. S. and Teicher H. (1978) {\it Probability Theory:
Independence, Interchangeability, Martingales.} New York: Springer.

\item{} Gell-Mann M. and Tsallis C. (2004) {\it Nonextensive Entropy -
Interdisciplinary Applications.} New York: Oxford University Press.

\item{} Gneiting T. (1997)  Normal scale mixtures and dual probability densities.
{\it J. Statist. Comput. Simul.} {\bf 59}, 375-384.

\item{} Gribbin D. W.,  Harris R. W. and Lau H. (1992)
Futures prices are not Stable-Paretian Distributed. {\it J of
Futures Markets.} {\bf 12}, 475-487.

\item{} Hahn M. G. and Zhang G. (1998)
Distinctions between the regular and empirical central limit
theorems for exchangeable random variables. In {\it Progress in
Probability}, {\bf 43}, Birkh\"auser Verlag Basel, Switzland,
111-143.

\item{} Hall J., Brorsen B. W. and Irwin S. (1989)
The distribution of futures prices: A test of the stable Paretian
and mixture of normals hypotheses. {\it J. of Financial and
Quantitative Analysis}. {\bf 24,} No 1, 105-116.

\item{} Hilhorst H. J. and Schehr G.
A note on q-Gaussians and non-Gaussians in statistical mechanics,
(2007) {\it J. Stat. Mech.} P06003.

\item{} Hanel R., Thurner S. and Tsallis C.  (2009) EPJ.
Scale-invariant correlated probabilistic model yields $q$-Gaussians
in the thermodynamic limit.

 \item {} Jiang X. and Hahn M. (2003)
 Central limit theorems for
exchangeable random variables when limits are mixture of normals.
{\it J. Theoret. Probab.} {\bf 16}, 543-571.

\item{} Kallenberg O. (2002) {\it  Foundations of Modern
Probability}, Springer Verlag.

\item{} Keilson J. and Steutel F. W. (1974)
Mixtures of distributions, moment inequalities and measures of
exponentiality and normality. {\it Ann. Probab.} {\bf 2,} 112-130.

\item{} Kon S. J. (1984)
Models of stock returns - A comparison. {\it J. Finance} {XXXIX}
147-165.

\item{} Marsh J. A.,  Fuentes M. A., Moyano L. G. and Tsallis C.
(2006) Influence of Global correlations on central limit theorems
and entropic extensivity. {\it Physica A}, {\bf 372,} 183-202.

\item{} Rodriguez A., Schwammle V. and Tsallis C. (2008) {\it Journal of
Statistical Mechanics}, P09006.

\item{} Tsallis C. (1988) Possible Generalization of
Boltzmann-Gibbs Statistics. {\it J. Statist. Phys.} {\bf 52},
479-487.

\item{} Tsallis C. 2009
{\it Introduction to nonextensive statistical mechanics.} Springer.

\item{} Umarov S, Tsallis C. and Steinberg S. (2008)
On a $q$-central limit theorem consistent with nonextensive
statistical mechanics. {\it Milan Journal of Mathematics,} {\bf 76}
307-328.

\item{} Umarov S. and Queiros S. M.- D. (2009)
Functional-differential equations for $F_q$ -transforms of
$q$-Gaussians. (ARXIV: 0802.0264)


\item{} Vignat C. and Plastino A. (2007)
Central limit theorem and deformed exponentials. {\it J. Phys. A:
Math. Theor.} {\bf 40}, F969-F978.


\ee


\end{document}